\documentstyle[12pt,amssymb]{article}

\def\fnote#1#2{\begingroup\def\thefootnote{#1}\footnote{#2}\addtocounter
{footnote}{-1}\endgroup}

\begin{document}

\baselineskip=18pt
\parskip=.01truein
\parindent=0pt

\hfill {\bf NSF-KITP-03-39}


\vskip 1truein

\centerline{{\large {\bf Black Hole Attractor Varieties and
Complex Multiplication}}\fnote{\diamond}{To appear in the
Proceedings of the {\it Workshop on Arithmetic, Geometry and
Physics around Calabi-Yau Varieties and Mirror Symmetry}, ed. by
J. Lewis and N. Yui, American Math. Society.}}

\vskip .5truein

\centerline{\sc Monika Lynker\fnote{\star}{Email:
mlynker@iusb.edu}$^1$, Vipul Periwal\fnote{\diamond}{Email:
vipul@gnsbiotech.com}$^2$ and Rolf
Schimmrigk\fnote{\dagger}{Email: netahu@yahoo.com,
rschimmr@kennesaw.edu}$^3$ }

\vskip .4truein

\centerline{{\it $^1$ Indiana University South Bend, South Bend,
IN 46634}}

\vskip .1truein

 \centerline{{\it $^2$ Gene Network Sciences, Ithaca, NY 14850}}

\vskip .1truein

\centerline{{\it $^3$ Kennesaw State University, Kennesaw, GA
30144}}

\vskip 0.9truein

\baselineskip=18pt \vskip .3truein

\centerline{\bf Abstract}

Black holes in string theory compactified on Calabi-Yau varieties
a priori might be expected to have moduli dependent features. For
example the entropy of the black hole might be expected to depend
on the complex structure of the manifold. This would be
inconsistent with known properties of black holes. Supersymmetric
black holes appear to evade this inconsistency by having moduli
fields that flow to fixed points in the moduli space that depend
only on the charges of the black hole.  Moore observed in the case
of compactifications with elliptic curve factors that these fixed
points are arithmetic, corresponding to curves with complex
multiplication. The main goal of this talk is to explore the
possibility of generalizing such a characterization to Calabi-Yau
varieties with finite fundamental groups.

\vskip .2truein

MSC~: Primary 14K22; Secondary 11R37


\renewcommand\thepage{}
\newpage

\baselineskip=22pt
\parskip=.1truein
\parindent=0pt
\pagenumbering{arabic}

\section{Introduction}

Arithmetic considerations have taken center stage in algebraic
geometry during the past two or three decades. In particular the
search for the still somewhat elusive concept of motives has
motivated much of current research. Few of these more modern
developments have had any impact on physics, even though it has
been clear for almost two decades that algebraic geometry is
central to the understanding of string theory. This is
 perhaps surprising because many of the number theoretic and
arithmetic results are closely linked to powerful analytic tools.
Recently however, arithmetic structures have been used to address
a variety of problems in string theory, such as the problem of
understanding aspects of the underlying conformal field theory of
Calabi-Yau manifolds \cite{s95-01}, the nature of black hole
attractor varieties \cite{m98}, and the behavior of periods under
reduction to finite fields \cite{cdv00}.

In this paper we describe some further developments and
generalizations of some of the observations made by Moore in his
analysis of the arithmetic nature of the so-called black hole
attractor varieties \cite{lps01}. The attractor mechanism
\cite{fks95, s96, fk96, fgk97} describes the radial evolution of
vector multiplet scalars of spherical dyonic black hole solutions
in $N=2$ supergravity coupled to abelian vector multiplets. Under
particular regularity conditions the vector scalars flow to a
fixed point in their target space. This fixed point is determined
by the charge of a black hole, described by a vector $\omega$ in
the lattice $\Lambda$
 of electric and magnetic charges of the $N=2$ abelian
gauge theory. If the $N=2$ supergravity theory is derived from a
type IIB string theory compactified on a Calabi-Yau space, the
vector multiplet moduli space is described by the moduli space
${\mathcal M}$ of complex structures of $X$, and the dyonic charge
vector takes values in the lattice $\Lambda = \rm H^3(X,\sf
Z\kern-.4em Z)$.

 One of the crucial observations made by
Moore is that in the context of simple toroidal product varieties,
such as the triple product of elliptic curves $E^3$, or the
product of the K3 surface and an elliptic curve K3$\times E$, the
attractor condition determines the complex moduli $\tau$ of the
tori to be given by algebraic numbers in a quadratic imaginary
field ${\mathbb Q}(\sqrt{D})$, $D<0$. This is of interest because
for particular points in the moduli space the elliptic curves
exhibit additional symmetries, that is, they admit so-called
complex multiplication (CM). For compactifications with such
toroidal factors Moore's analysis then appears to indicate a
strong link between the 'attractiveness' of varieties in string
theory and their complex multiplication properties. We will
briefly review Moore's observations in Section 2.

Calabi-Yau varieties with elliptic factors are very special
because they have infinite fundamental group, a property not
shared by Calabi-Yau manifolds in general. Other special features
of elliptic curves are not present in general either. In
particular Calabi-Yau spaces are not abelian varieties and they do
not, in any obvious fashion, admit complex multiplication. Hence
it is not clear how Moore's observations can be generalized. It is
this problem which is addressed in \cite{lps01}. In order to
formulate such a generalization we adopt a cohomological approach
and view the modular parameter of the elliptic curve as part of
the primitive cohomology. In the case of elliptic curves $E$ this
is simply a choice of point of view because there exists an
isomorphism between the curve itself and its Jacobian defined by
$J(E)={\rm H}^1(E,{\mathbb C})/{\rm H}^1(E,{\mathbb Z})$ described
by the Abel-Jacobi map $j: E \rightarrow J(E)$. These varieties
are abelian.

The Jacobian variety of elliptic (and more general) curves has a
natural generalization to higher dimensional varieties defined by
the intermediate Jacobian of Griffiths. It would be natural to use
Griffiths' construction in an attempt to generalize the elliptic
results described above. In general, however, the intermediate
Jacobian is not an abelian variety and does not admit complex
multiplication. For this reason we will proceed differently by
constructing a decomposition of the intermediate cohomology of the
Calabi-Yau variety. We then use this decomposition to formulate a
generalization of the concept of complex multiplication of black
hole attractor varieties. To achieve this we formulate complex
multiplication in this more general context by analyzing in some
detail the cohomology group $\rm H^3(X)$  of weighted Fermat
hypersurfaces $X$.

The paper is organized as follows. In Section 4 we briefly review
the necessary background of abelian varieties. In Section 5 we
show how abelian varieties can be derived from Calabi-Yau
hypersurfaces by showing that the cohomology of such varieties can
be constructed from the cohomology of curves embedded in these
higher dimensional varieties. This leads us to abelian varieties
defined by the Jacobians of curves.
 Such abelian varieties do not, in general, admit complex
multiplication. However, it is known that Jacobians of ordinary
projective Fermat curves split into abelian factors which do admit
complex multiplication. We briefly describe this construction and
demonstrate that this property generalizes to Brieskorn-Pham
curves. Combining these results shows that we can consider the
complex multiplication type of Calabi-Yau varieties as determined
by the CM type of their underlying Jacobians.

When these results are applied to describe the complex
multiplication type of a particularly simple black hole attractor
variety it emerges that its complex multiplication leads precisely
to the field determined by its periods. It is in fact not
completely unexpected that we might be able to recover the field
of periods by considering the complex multiplication type. The
reason for this is a conjecture of Deligne \cite{d79} which states
that the field determined by the periods of a critical motive is
determined by its L-function. Because Deligne's conjecture is
important for our general view of the issue at hand we briefly
describe this conjecture in  Section 3 in order to provide a
broader perspective. It is important to note that Deligne's
conjecture is in fact a theorem in the context of pure projective
Fermat hypersurfaces \cite{b86}, but has not been proven in the
context of weighted hypersurfaces. Our results in essence can be
viewed as support of this conjecture even in this more general
context. In Section 2 we briefly review the physical setting of
black hole attractors in type IIB theories, as well as Moore's
solution of the ${\rm K3} \times E$ solution of the attractor
equations. In Section 7 we summarize our results and indicate
possible generalizations.

\section{Attractor Varieties}

\subsection{Compactified type IIB string theory}

We consider type IIB string theory compactified on Calabi-Yau
threefold varieties. The field content of string theory in 10D
space $X_{10}$ splits into two sectors according to the boundary
conditions on the world sheet. The Neveu-Schwarz fields are given
by the metric $g \in \Gamma(X_{10}, \rm T^* X_{10}\otimes \rm T^*
X_{10})$, an antisymmetric tensor field $B \in \Gamma(X_{10},
\Omega^2))$ and the dilaton scalar $\phi \in C^{\infty}(X_{10},
{\mathbb R})$. The Ramond sector is spanned by even antisymmetric
forms $A^p \in \Gamma(X_{10}, \Omega^p)$ of rank $p$ zero, two,
and four. Here $\Omega^p \longrightarrow X$ denotes the bundle of
$p$-forms over the variety $X$.

In the context of the black hole solutions considered in
\cite{fks95} the pertinent sectors are given by the metric and the
five-form field strength ${\bf F}$ of the Ramond-Ramond 4-form
$A^4$. The metric is assumed to be static, spherically symmetric,
asymptotically Minkowskian, and should describe extremally charged
black holes, leading to the ansatz \begin{equation} ds^2
=-e^{2U(r)} dt\otimes dt + e^{-2U(r)}(dr\otimes dr +
r^2\sigma_2),\end{equation} where $r$ is the spatial three
dimensional radius, $\sigma_2$ is the 2D angular element, and the
asymptotic behavior is encoded via $e^{-U(r)} \rightarrow \infty$
for $r \rightarrow \infty$. The ten-dimensional five-form
$\mathbf{F}$ leads to a number of different four-dimensional
fields, the most important in the present context being the field
strengths $F^L$ of the four dimensional abelian fields, the number
of which depends on the dimension of the cohomology group $\rm
H^3(X)$ via $A^4_{\mu mnp}(x,y) = \sum_L A^{4L}_{\mu}(x)
\omega^L_{mnp}(y)$, where $\{\omega_L\}_{L=1,...,b_3}$ is a basis
of $\rm H^3(X)$. This is usually written in a symplectic basis
$\{\alpha_a, \beta^a\}_{a=0,...,h^{2,1}}$, for which $\int_X
\alpha_a \wedge \beta^b =\delta_a^b$, as an expansion of the field
strength
\begin{equation} {\bf F} (x,y) = F^a(x) \wedge \alpha_a -
G_a(x) \wedge
\beta^a. \end{equation} Being a five-form in ten dimensions, the
field strength ${\bf F}$ admits (anti)self-duality constraints
with respect to Hodge duality, ${\bf F} = \pm \ast_{10} {\bf F}$.
The ten dimensional Hodge operator $\ast_{10}$ factorizes into a
4D and a 6D part $\ast_{10} = \ast_4 \ast_6$. A solution to the
antiselfduality constraint in 10D as well as the Biachi identity
$d\mathbf{F} =0$ can be obtained by setting \cite{m98}
\begin{equation} \mathbf{F} = \textrm{Re} \left({\bf E}\wedge
(\omega^{2,1} + \omega^{0,3}
) \right), \end{equation} where \cite{f97,d98} \begin{equation}
{\bf E} \equiv q \sin \theta d\theta \wedge d\phi - i q
\frac{e^{2U(r)}}{r^2} dt \wedge dr \end{equation} is a 2-form for
which the four-dimensional Hodge duality operator leads to
$\ast_{4} {\bf E} =i{\bf E}$.

\subsection{}
The dynamics of a string background configuration can be derived
by either reducing the IIB effective action with a small
superspace ansatz \cite{fgk97}, or via the supersymmetry variation
constraints of the fermions in nontrivial backgrounds, in
particular the gravitino and gaugino variations. Defining on $\rm
H^3(X)$ an inner product $<\cdot , \cdot >$ via
\begin{equation} <\omega, \eta> = \int_X \omega \wedge \eta,
\end{equation} the gravitino equation involves the integrated
version of the 5-form field strength \cite{bcdffrsv95}
\cite{cdf95}
\begin{equation} {\rm T}^{-} = e^{K/2} <\Omega, {\bf {\rm F}}^{-}>
=  e^{K/2} \left( {\mathcal G}_a F^{-,a}(x) - z^a G_a^{-}(x)
\right)
\end{equation} with the K\"ahler potential \begin{equation}
e^{-K} = i<\Omega, \bar{\Omega}> = -i(z^a \bar{{\mathcal G}}_a -
\bar{z}^a {\mathcal G}_a) ,\end{equation} where the second
equation is written in terms of the periods $z^a = <\Omega,
\beta^a> = \int_{A^a}\Omega$, and ${\mathcal G}_a = <\Omega,
\alpha_a> = \int_{B_a} \Omega$ with respect to a symplectic dual
homological basis $\{A^a, B_a\}$, whose dual cohomological basis
is denoted by $\{\alpha_a, \beta^a\} \subset \rm H^3(X)$. The
holomorphic three-form thus can be expanded as $\Omega =
z^a\alpha_a -{\mathcal G}_a \beta^a$.

The supersymmetry transformation of the gravitino $\psi^A =
\psi_{\mu}^A dx^{\mu}$ can then be written as \begin{equation}
\delta \psi^A = D \varepsilon^A + dx^{\mu} \mathrm{T}^{-}_{\mu
\nu} \gamma^{\nu} (\epsilon \varepsilon)^A,\end{equation} where
$\gamma^{\mu}$  denote the covariant Dirac matrices.  The
variation of the gaugino of the abelian multiplets takes the form
\begin{equation} \delta \lambda^{iA} = i\gamma^{\mu}
\partial_{\mu} z^i\epsilon^A +
\frac{i}{2}G^{-,i}_{\mu \nu} \gamma^{\mu \nu} (\epsilon
\varepsilon )^A.\end{equation}

\subsection{}
Plugging these ingredients into the supersymmetry transformation
behavior of the gravitino and the gaugino fields, and demanding
that the vacuum remains fermion free, leads to  the following
equations for the moduli and the spacetime
function $U(r)$ \begin{eqnarray} \frac{dU}{d\rho} &=& - e^U|Z|
\nonumber \\
\frac{dz^i}{d\rho} &=& -2e^U g^{i\bar{j}} \partial_{\bar{j}}
|Z|,\end{eqnarray} where
\begin{equation} Z(\Gamma) = e^{K/2} \int_{\Gamma} \Omega
=  e^{K/2} \int_X
\eta_{\Gamma} \wedge \Omega \end{equation} is the central charge
of the cycle $\Gamma \in \rm H_3(X)$ with Poincare dual
$\eta_{\Gamma} \in \rm H^3(X)$ and $g_{i\bar{j}} = \partial_i
\partial_{\bar{j}} K$ is the metric derived from the K\"ahler
potential $K$.

The fixed point condition of the attractor equation can be
rewritten in a geometrical way as the Hodge condition
\begin{equation} \rm H^3(X,{\mathbb Z}) \ni \omega =
\omega^{3,0} + \omega^{0,3}. \end{equation} Writing $\omega^{3,0}
=-i\bar{C} \Omega$ this can also be formulated as
\begin{eqnarray} ip^a &=& \bar{C} z^a - C \bar{z}^a \nonumber \\
     iq_a &=& \bar{C} {\mathcal G}_a - C\bar{{\mathcal G}}_a,
     \end{eqnarray}
where $C=e^{K/2}Z$. This system describes a set of $b_3(X)$
charges $(p^a,q_a)$ determined by the physical 4-dimensional input
which in turn determines the system of complex periods of the
Calabi-Yau variety. Hence the equations should be solvable. The
interesting structure of the fixed point which emerges is that the
central charges are determined completely in terms of the charges
of the four-dimensional theory. As a consequence the 4D geometry
is such that the horizon is a moduli independent quantity. This is
precisely as expected because the black hole entropy should not
depend on adiabatic changes of the environment \cite{lw95}.

\subsection{}
In reference \cite{m98} it is noted that two types of solutions of
the attractor equations have particularly interesting properties.
The first of these is provided by the triple product of a torus,
while the second is a product of a K3 surface and a torus. Both
solutions are special in the sense that they involve elliptic
curves. In the case of the product threefold $X=\mathrm{K3} \times
E$ the simplifying feature is that via K\"unneth's theorem one
finds ${\rm H^3}(\mathrm{K3} \times E) \cong {\rm H}^2(K3) \otimes
{\rm H}^1(E)$, and therefore the cohomology group of the threefold
in the middle dimension is isomorphic to two copies of the
cohomology group $\rm H^2(\mathrm{K3})$. The attractor equations
for such threefolds have been considered in \cite{adf96}. The
resulting constraints determine the holomorphic form of both
factors in terms of the charges $(p,q)$ of the fields. Moore finds
that the complex structure $\tau$ of the elliptic curve $E =
{\mathbb C}/{\mathbb Z} \tau +{\mathbb Z}$
 is solved as
\begin{equation} \tau_{p,q} = \frac{p\cdot q +\sqrt{D_{p,q}}}{p^2},
\end{equation} where
$D_{p,q}=(p\cdot q)^2 - p^2q^2$ is the discriminant of a BPS state
labelled by \begin{equation} \omega = (p,q) \in {\rm H}^3
(\mathrm{K3}\times E,{\mathbb Z}).\end{equation} The holomorphic
two form on K3 is determined as $\Omega^{2,0} = {\mathcal C}
(q-\bar{\tau} p)$, where ${\mathcal C}$ is a constant. Moore makes
the interesting observation that this result is known to imply
that the elliptic curve determined by the attractor equation is
distinguished by exhibiting a particularly symmetric structure,
i.e. that the endomorphism algebra $\mathrm{End}(E)$ is enlarged.
In general $\mathrm{End}(E)$ is just the ring ${\mathbb Z}$ of
rational integers. For special curves however there are two other
possibilities for which $\mathrm{End}(E)$ is either an order of a
quadratic imaginary field, or it is a maximal order in a
quaternion algebra. The latter possibility can occur only when the
field $K$ over which $E$ is defined has positive characteristic.
Elliptic curves are said to admit complex multiplication if the
endomorphism algebra is strictly larger than the ring of rational
integers.

\subsection{}
The important point here is that the property of complex
multiplication appears if and only if the $j$-invariant $j(\tau)$
is an algebraic integer. This happens if and only if the modulus
$\tau$ is an imaginary quadratic number. The $j$-invariant of the
elliptic curve $E_{\tau}$ can be defined in terms of the
Eisenstein series \begin{equation} E_k(\tau) = \frac{1}{2}
\sum_{\stackrel{m,n \in {\mathbb Z}}{m,n~{\rm coprime}}}
\frac{1}{(m\tau +n)^k}
\end{equation} as \begin{equation} j(\tau) =
\frac{E_4(\tau)^3}{\Delta(\tau)},\end{equation}
 where
$1728 \Delta(\tau) = E_4(\tau)^3 - E_6(\tau)^2$. In general the
$j$-function does not take algebraic values, not to mention values
in an imaginary quadratic field. We therefore see that in this
elliptic setting the solutions of the attractor equations can be
characterized as varieties which admit complex multiplication and
determine a quadratic imaginary field $K_D = {\mathbb
Q}(i\sqrt{|D|})$.

Once this is recognized several classical results about elliptic
curves with complex multiplication are available to illuminate the
nature of the attractor variety. One of these results is that the
extension $K_D(j(\tau))$ obtained by adjoining the $j$-value to
$K_D$ is the Hilbert class field. Geometrically there is a
Weierstrass model, i.e. a projective embedding of the elliptic
curve of the form \begin{equation} y^2 + a_1xy +a_3y = x^3 +
a_2x^2 + a_4x + a_6 \end{equation} that is defined over this
extension $K_D(j(\tau))$.

Even more interesting is that it is possible to construct from the
geometry of the elliptic curve the maximal abelian extension of
$K_D$ by considering the torsion points $E_{\mathrm{tor}}$ on the
curve $E$, i.e. points of finite order with respect to the group
law. To do this consider the Weber function $\phi_{E}$ on the
curve $E$. Assuming that the characteristic of the field $K_D$ is
different from 2 or 3, the elliptic curve can be embedded via the
simplified Weierstrass form \begin{equation} y^2=x^3 +Ax +B
\end{equation} with discriminant
\begin{equation} \Delta = -16(4A^3+27B^2) \neq 0. \end{equation}
The Weber function can then be defined as \begin{equation}
\phi_{E}(p) = \left\{
 \begin{tabular}{l l}
     $\frac{AB}{\Delta}x(p)$  & if $j(E) \neq 0 ~
     \mathrm{or}~ 1728$ \\
     $\frac{A^2}{\Delta}x^2(p)$  & if $j(E)=1728$  \\
     $\frac{B}{\Delta}x^3(p)$   & if $j(E)=0$.  \\
  \end{tabular} \right\}
\end{equation} and the Hilbert class field
$K_D(j(\tau))$ can be extended to the full maximal abelian
extension $K_D^{\rm ab}$ of $K_D$ by adjoining the Weber values of
the torsion points $K_D^{\rm ab} = K_D(j(\tau), \{\phi_{E}(T)~|~T
\in \mathrm{E}_{\mathrm{tor}}\})$. We see from this that the
attractor equations pick out special elliptic curves which lead to
a rich arithmetic structure. It is this set of tools which we wish
to generalize to the framework of Calabi-Yau varieties proper,
i.e. those with finite fundamental group.

\section{Deligne's Period Conjecture}

\subsection{}
In this section we briefly review Deligne's conjecture in its
motivic formulation. This is useful because it will allow us to
provide a general perspective for our results which will furnish a
useful general framework in which to investigate the arithmetic
nature  of attractor varieties. Motives are somewhat complicated
objects whose status is reminiscent to string theory: different
realizations are used to probe what is believed to be some yet
unknown unifying universal cohomology theory of varieties. More
precisely, motives are characterized by a triplet of different
cohomology theories together with a pair of compatibility
homomorphisms. In terms of these ingredients a motive then is
described by the quintuplet of objects
\begin{equation} (\mathrm{M}_{\mathrm{B}}, \mathrm{M}_{\mathrm{dR}},
\mathrm{M}_{\ell}, \mathrm{I}_{\mathrm{B}, \sigma},
\mathrm{I}_{\ell, \bar{\sigma}}),\end{equation} where the three
first entries are cohomology objects constructed via Tate twists
from the Betti, de Rham, and \'etale cohomology, respectively.
Furthermore $\mathrm{I}_{\mathrm{B},\sigma}$ describes a map
between the Betti and de Rham cohomology, while $\mathrm{I}_{\ell,
\bar{\sigma}}$ is a map between Betti and \'etale
cohomology\footnote[2]{Detailed reviews of motives can be found in
\cite{jks94}.}. In the following we will focus mostly on motives
derived from the first (co)homology groups $\mathrm{H}^1(A)$ and
$\mathrm{H}_1(A)$ of abelian varieties $A$, as well as the
primitive cohomology of Fermat hypersurfaces.

\subsection{}
The second ingredient in Deligne's conjecture is the concept of a
geometric $L-$function. This can be described in a number of
equivalent ways. Conceptually the perhaps simplest approach
results when it can be derived via Artin's zeta function as the
Hasse-Weil L-function induced by the underlying variety, i.e. a
way of counting solutions of the variety over finite fields. The
complete L-function receives contributions from two fundamentally
different factors $\Lambda(\mathrm{M},s) =
L_{\infty}(\mathrm{M},s) L(\mathrm{M},s)$. The infinity term
$L_{\infty}(M,s)$ originates from those fields over which the
underlying variety has bad reduction, i.e. it is singular, while
the second term $L(\rm M,s)$ collects all the information obtained
from the finite fields over which the variety is smooth. The
complete L-function is in general expected to satisfy a functional
equation, relating its values at $s$ and $1-s$. A motive is called
critical if neither of the infinity factors in the functional
equation has a pole at $s=0$.

\subsection{} The final ingredient is the concept of the period of a
motive, a generalization  of ordinary periods of varieties.
Viewing the motive $M$ as a generalized cohomology theory Deligne
formulates the notion of a period $c^+(M) \in {\mathbb
C}^*/{\mathbb Q}^*$ by taking the determinant of a compatibility
homomorphism
\begin{equation} I_{\mathrm{B},\sigma}: \mathrm{M}_{\mathrm{B}}
\longrightarrow \mathrm{M}_{\mathrm{dR}}
\end{equation} between the Betti and the deRham realizations of
the motive $M$. Deligne's basic conjecture then relates the period
and the L-function  via $L(M,0)/c^+(M) \in {\mathbb Q}$. Contact
with the Hasse-Weil L-function is made by noting that for motives
of the type $M=\mathrm{H}(X)(m)$ with Tate twists one has
$L(\mathrm{M},0)=L(X,m)$.

\subsection{}
Important for us is a generalization of this conjecture that
involves motives with coefficients in a field $E$. Such motives
can best be described via algebraic Hecke characters which are of
particular interest for us because they come up in the L-function
of projective Fermat varieties. Algebraic Hecke characters were
first introduced by Weil as Hecke characters of type $A_0$, which
is what they are called in the older literature. In the context of
motives constructed from these characters the field $E$ becomes
the field of complex multiplication.  In this more general context
Deligne's conjecture says that the L-function of the motive and
the period take values in the same field, the CM field of the
motive.

{\bf Deligne Period Conjecture:} \begin{equation}
\frac{L(M,0)}{c^+(M)} \in E.\end{equation}

For Fermat hypersurfaces Deligne's conjecture is in fact a
theorem, proven by Blasius \cite{b86}.

\section{Abelian Varieties from Brieskorn-Pham Hypersurfaces}

\subsection{}
An abelian variety over some number field $K$ is a smooth,
geometrically connected, projective variety which is also an
algebraic group with the group law $A\times A \longrightarrow A$
defined over $K$. A concrete way to construct abelian varieties is
via complex tori ${\mathbb C}^n/\Lambda$ with respect some lattice
$\Lambda$ that is not necessarily integral and admits a Riemann
form. The latter is defined as an ${\mathbb R}$-bilinear form
$<,>$ on ${\mathbb C}^n$ such that $<x,y>$ takes integral values
for all $x,y \in \Lambda$, $<x,y>=-<y,x>$, and $<x,iy>$ is a
positive symmetric form, not necessarily non-degenerate. Then one
has the result that a complex torus ${\mathbb C}^n/\Lambda$ has
the structure of an abelian variety if and only if there exists a
non-degenerate Riemann form on ${\mathbb C}^n/\Lambda$.

\subsection{}
A special class of abelian varieties are those of CM-type,
so-called complex multiplication type. The reason for these
varieties to be of particular interest is that certain number
theoretic question can be addressed in a systematic fashion for
this class. Consider a number field $K$ over the rational numbers
${\mathbb Q}$ and denote by $[K:{\mathbb Q}]$ the degree of the
field $K$ over ${\mathbb Q}$, i.e. the dimension of $K$ over the
subfield ${\mathbb Q}$. An abelian variety $A$ of dimension $n$ is
called a CM$-$variety if there exists an algebraic number field
$K$ of degree $[K:{\mathbb Q}]=2n$ over the rationals ${\mathbb
Q}$ which can be embedded into the endomorphism algebra
$\mathrm{End}(A) \otimes {\mathbb Q}$ of the variety. More
precisely, a CM-variety is a pair $(A,\theta)$ with $\theta: K
\longrightarrow \mathrm{End}(A) \otimes {\mathbb Q}$ an embedding
of $K$. It follows from this that the field $K$ necessarily is a
CM field, i.e. a totally imaginary quadratic extension of a
totally real field.  The important ingredient here is that the
restriction to $\theta(K) \subset \mathrm{End}(A)\otimes {\mathbb
Q}$ is equivalent to the direct sum of $n$ isomorphisms
$\varphi_1,...,\varphi_n \in \mathrm{Iso}(K,{\mathbb C})$ such
that $\mathrm{Iso}(K,{\mathbb C}) = \{\varphi_1,...,\varphi_n,
\rho\varphi_1,....,\rho \varphi_n\}$, where $\rho$ denotes complex
conjugation. These considerations lead to the definition  of
calling the pair $(K,\{\varphi_i\})$ a CM-type, in the present
context, the CM-type of a CM-variety $(A,\theta)$.

\subsection{} The context in which these concepts will appear below
is provided by varieties which have complex multiplication by a
cyclotomic field $K={\mathbb Q}(\mu_n)$, where $\mu_n$ denotes the
cyclic group generated by a nontrivial $n$'th root of unity
$\xi_n$. The degree of ${\mathbb Q}(\mu_n)$ is given by $[{\mathbb
Q}(\mu_n):{\mathbb Q}]=\phi(n)$, where $\phi(n)=\# \{m\in {\mathbb
N}~|~m<n, ~\mathrm{gcd}(m,n)=1\}$ is the Euler function. Hence the
abelian varieties encountered below will have complex dimension
$\phi(n)/2$.

In the following we first reduce the cohomology of the
Brieskorn-Pham varieties to that generated by curves and then
analyze the structure of the resulting weighted curve Jacobians.

\subsection{Curves and the cohomology of threefolds}

The difficulty of general higher dimensional varieties is that
there is no immediate way to recover abelian varieties, thus
making it non-obvious how to generalize the concept of complex
multiplication from elliptic curves with complex multiplication,
which are abelian varieties. As a first step we need to
disentangle the Jacobian of the elliptic curve from the curve
itself. This would lead us to the concept of the
middle-dimensional cohomology, more precisely the intermediate
(Griffiths) Jacobian which is the appropriate generalization of
the Jacobian of complex curves. The problem with this intermediate
Jacobian is that it is not, in general, an abelian variety.

We will show now that it is possible nevertheless to recover
abelian varieties as the basic building blocks of the intermediate
cohomology in the case of weighted projective hypersurfaces. The
basic reason for this is that the cohomology $\rm H^3(X)$ for
these varieties decomposes into the monomial part and the part
coming from the resolution. The monomial part of the intermediate
cohomology can easily be obtained from the cohomology of a
projective hypersurface of the same degree by realizing the
weighted projective space as a quotient variety with respect to a
product of discrete groups determined by the weights of the
coordinates. For projective varieties ${\mathbb P}_n[d]$ it was
shown in \cite{sk79} that the intermediate cohomology can be
determined by lower-dimensional varieties in combination with Tate
twists. Denote the Tate twist by \begin{equation} \rm
H^i(X)(j):=\rm H^i(X)\otimes W^{\otimes j}\end{equation} with
$W=\rm H^2({\mathbb P}_1)$ and let $X_d^{r+s}$ be a Fermat variety
of degree $d$ and dimension $r+s$. Then \begin{eqnarray} \rm
H^{r+s}(X_d^{r+s}) & \oplus &\sum_{j=1}^r \rm
H^{r+s-2j}(X_d^{r-1})(j) \oplus \sum_{k=1}^s
\rm H^{r+s-2k}(X_d^{s-1})(k) \nonumber \\
 & & ~~~~ \cong \rm H^{r+s}(X_d^r\times
X_d^s)^{\mu_d} \oplus \rm H^{r+s-2}(X_d^{r-1}\times
X_d^{s-1})(1).\end{eqnarray} Here $\mu_d$ is the cyclic group of
order $d$ which acts on the individual factors as \begin{equation}
[(x_0,...,x_r),(y_0,...,y_s)] \mapsto [(x_0,...,x_{r-1},\xi x_r),
(y_0,...,y_{s-1},\xi y_s)] \end{equation} and induces an action on
the hypersurfaces.

\subsection{} This leaves only the part of the intermediate
cohomology of the weighted hypersurface that originates from the
resolution. It was shown in \cite{cls90} that the only singular
sets on arbitrary weighted hypersurface Calabi-Yau threefolds are
either points or curves. The resolution of singular points
contributes to the even cohomology group $\rm H^2(X)$ of the
variety, but does not contribute to the middle-dimensional
cohomology group $H^3(X)$. Hence we need to be concerned only with
the resolution of curves. This can be described for general CY
hypersurface threefolds as follows. If a discrete symmetry group
${\mathbb Z}/n{\mathbb Z}$ of order $n$ acting on the threefold
leaves a curve invariant then the normal bundle has fibres
${\mathbb C}_2$ and the discrete group induces an action on these
fibres which can be described by a matrix
\begin{displaymath}
\left( \begin{array}{cc} \alpha^{mq} & 0 \\
0 & \alpha^m\\
\end{array}\right),
\end{displaymath}
where $\alpha$ is an $n$'th root of unity and $(q,n)$ have no
common divisor. The quotient ${\mathbb C}_2/({\mathbb Z}/n{\mathbb
Z})$ by this action has an isolated singularity which can be
described as the singular set of the surface in ${\mathbb C}_3$
given by the equation
\begin{equation} S=\{(z_1,z_2,z_3)\in {\mathbb C}_3 ~|~
z_3^n=z_1z_2^{n-q}\}.\end{equation} The resolution
of such a singularity is completely determined by the type $(n,q)$
of the action by computing the continued fraction of $\frac{n}{q}$
\begin{equation} \frac{n}{q}= b_1 - \frac{1}{b_2 - \frac{1}{\ddots
- \frac{1}{b_s}}} \equiv [b_1,...,b_s].\end{equation} The numbers
$b_i$ specify completely the plumbing process that replaces the
singularity and in particular determine the additional generator
to the cohomology $\rm H^*(X)$ because the number of ${\mathbb
P}_1$s introduced in this process is precisely the number of steps
needed in the evaluation of $\frac{n}{q}=[b_1,...,b_s]$. This can
be traced to the fact that the singularity is resolved by a bundle
which is constructed out of $s+1$ patches with $s$ transition
functions that are specified by the numbers $b_i$. Each of these
glueing steps introduces a sphere, which in turn supports a
(1,1)-form. The intersection properties of these 2-spheres are
described by Hirzebruch-Jung trees, which for a ${\mathbb
Z}/n{\mathbb Z}$ action is just an SU($n+1$) Dynkin diagram, while
the numbers $b_i$ describe the intersection numbers. We see from
this that the resolution of a curve of genus $g$ thus introduces
$s$ additional generators to the second cohomology group $\rm
H^2(X)$, and $g\times s$ generators to the intermediate cohomology
$\rm H^3(X)$.

Hence we have shown that the cohomology of weighted hypersurfaces
is determined completely by the cohomology of curves. Since the
Jacobian, which we will describe in the next subsection, is the
only motivic invariant of a smooth projective curve this says that
for weighted hypersurfaces the main motivic structure is carried
by their embedded curves.

\subsection{Cohomology of weighted curves}

For smooth algebraic curves $C$ of genus $g$ the de Rham
cohomology group $\rm H^1_{\rm dR}(C)$ decomposes (over the
complex number ${\mathbb C}$) as \begin{equation} \rm H^1_{\rm
dR}(C)~\cong~\rm H^0(C, \Omega^1) \oplus \rm H^1(C,{\mathcal
O}).\end{equation} The Jacobian $J(C)$ of a curve $C$ of genus $g$
can be identified with $J(C)={\mathbb C}^g/\Lambda$, where
$\Lambda$ is the period lattice
\begin{equation} \Lambda:= \left\{(\dots,\int_a \omega_i,\dots)~|~
a \in \rm H_1(C,{\mathbb Z}),~ \omega_i \in \rm H^1(C)
\right\},\end{equation} where the $\omega_i$ form a basis. Given a
fixed point $p\in C$ on the curve there is a canonical map from
the curve to the Jacobian defined as
\begin{equation} \varphi_0: C \longrightarrow J(C) \end{equation}
via \begin{equation} p \mapsto \left(\dots, \int_{p_0}^p
\omega_{r,s},\dots\right)~ \mathrm{mod}~ \Lambda .\end{equation}

We are interested in curves of Brieskorn-Pham type, i.e. curves of
the form \begin{equation} {\mathbb P}_{(1,k,\ell)}[d]\ni
\left\{x^d + y^a + z^b =0 \right\},\end{equation} such that
$a=d/k$ and $b=d/\ell$ are positive rational integers. Without
loss of generality we can assume that $(k,\ell)=1$. We claim that
for smooth elements in ${\mathbb P}_{(1,k,\ell)}[d]$ the set of
forms
 \begin{eqnarray}
 & & \Omega({\mathbb P}_{(1,k,\ell)}[d])   \nonumber \\
 & & ~~~~ = \left\{ \omega_{rst}= y^{s-1}z^{t-d/\ell}dy~
 {\Large |}~ r+ks+\ell t = 0~{\mathrm{mod}}~d, ~
 \left(\begin{array}{c}1\leq r \leq d-1, \\
~1\leq s \leq \frac{d}{k}-1,\\ ~1\leq t \leq \frac{d}{\ell} -1 \\
\end{array}\right)\right\} \nonumber \\
\label{weighted-basis}\end{eqnarray}
  defines a basis
for the de Rham cohomology group $\rm H^1_{\rm dR}({\mathbb
P}_{(1,k,\ell)}[d])$.

In order to show this we view the projective space as the quotient
with respect to the actions ${\mathbb Z}_k:[0~1~0]$ and ${\mathbb
Z}_{\ell}: [0~0~1]$. This allows us to view the weighted curve as
the quotient of a pure projective Fermat type curve
\begin{equation} {\mathbb P}_{(1,k,\ell)}[d] =
{\mathbb P}_2[d]/{\mathbb Z}_k \times
{\mathbb Z}_{\ell}: \left[\begin{array}{ccc} 0&1&0\\
0&0&1\\
\end{array}\right].\end{equation} These weighted curves are smooth
and hence their cohomology is determined by considering those
forms on the pure projective curve ${\mathbb P}_2[d]$ which are
invariant with respect to the group actions. A basis for
$\Omega({\mathbb P}_2[d])$ is given by the set of forms
\begin{equation} \left\{ \omega_{rst}= y^{s-1}z^{t-d}dy~{\Large |}~
0<r,s,t<d,~~r+s+t=0~(\mathrm{mod} ~d),~~r,s,t\in {\mathbb
N}\right\}.\label{basis2} \end{equation}Denote the generator of
the ${\mathbb Z}_{\ell}$ action by $\alpha$ and consider the
induced action on $\omega_{rst}$ \begin{equation} \omega_{rst}
\mapsto \alpha^s \omega_{rst}.\end{equation} It follows that the
only forms that descend to the quotient with respect to ${\mathbb
Z}_{\ell}$ are those for which $r = 0 (\mathrm{mod}~\ell)$.
Similarly we denote by $\beta$ the generator of the action
${\mathbb Z}_k$ and consider the induced action on the forms
$\omega_{rst}$
\begin{equation} \omega_{rst} \mapsto \beta^{t-d} \omega_{rst}.
\end{equation}
Again we see that the only forms that descend to the quotient are
those for which $s=0(\mathrm{mod}~k)$.

\subsection{Abelian Varieties of weighted Jacobians}

It was shown by Faddeev  \cite{f61}\footnote[4]{More accessible
references of the subject are \cite{w76} \cite{g78}.}
 in the case of Fermat curves
that the Jacobian $J(C_d)$ splits into a product of abelian
factors \begin{equation} J(C_d) \cong \prod_{{\mathcal O}_i \in
I/{\mathbb Z}_d^*} A_{{\mathcal O}_i}, \end{equation} where the
${\mathcal O}_i$ are orbits in $I$ of the multiplicative subgroup
${\mathbb Z}_d^*$ of the group ${\mathbb Z}_d \equiv {\mathbb
Z}/d{\mathbb Z}$. More precisely it was shown that there is an
isogeny
\begin{equation} i: J(C_d) \mapsto \prod_{{\mathcal O}_i
\in I/{\mathbb Z}_d^*} A_{{\mathcal O}_i},\end{equation} where an
isogeny $i: A \rightarrow B$ between abelian varieties is defined
to be a surjective
 homomorphism with finite kernel.
We adapt this discussion to the weighted case.

 Denote the index set of triples $(r,s,t)$
parametrizing all one-forms  by ${\mathcal I}$. The cyclic group
${\mathbb Z}_d$ again acts on ${\mathcal I}$ and the
multiplicative subgroup $({\mathbb Z}_d)^*$ produces a set of
orbits
\begin{equation} [(r,s,t)] \in {\mathcal I}/({\mathbb Z}_d)^*.
\end{equation} Each of these orbits
leads to an abelian variety $A_{[(r,s,t)]}$ of dimension
\begin{equation} \mathrm{dim} A_{[(r,s,t)]} = \frac{1}{2}
\phi\left(\frac{d}{\mathrm{gcd}(r,ks,\ell t,d)} \right) ,
\end{equation} and complex multiplication with respect to
the field \begin{equation} K_{[(r,s,t)]} = {\mathbb
Q}(\mu_{d/\mathrm{gcd}(r,ks,\ell t,d)}).\end{equation} This leads
to an isogeny
\begin{equation}  i: J({\mathbb P}_{(1,k,\ell)}[d]) \mapsto
\prod_{[(r,s,t)]\in  {\mathcal I}/{\mathbb Z}_d^*}
A_{[(r,s,t)]}.\end{equation}

 The complex
multiplication type of the abelian factors $A_{[(r,s,t)]}$ of the
Jacobian $J(C)$ can be identified with the set \begin{equation}
\rm H_{rst} := \left\{a\in ({\mathbb Z}/d{\mathbb
Z})^*~|~<ar>+<aks>+<a\ell t>=d\right\}
\end{equation} via a homorphism from $\rm H_{rst}$ to the Galois
group. This leads to the subgroup of the Galois group of the
cyclotomic field given by \begin{equation} \mathrm{Gal}({\mathbb
Q}(\mu_{rst})/{\mathbb Q}) \supset {\mathcal H}_{rst} =
\left\{\sigma \in \mathrm{Gal}({\mathbb Q}(\mu_{rst})/{\mathbb
Q})~|~a \in \mathrm{H}_{rst}\right\}
\end{equation} and the CM type of $(A_{[(r,s,t)]},\theta_{rst})$
can be given as
\begin{equation} \left({\mathbb Q}(\mu_{rst}), \{\varphi_1,\cdots,
\varphi_n\} = {\mathcal H}_{rst}\right),\end{equation} where
$n=\phi(d_{rst})/2$.

\section{Summary and Generalizations}

We have seen that the concepts used to describe attractor
varieties in the context of elliptic compactifications can be
generalized to Calabi-Yau varieties with finite fundamental
groups. We have mentioned above that the abelian property is
neither carried by the variety itself nor the generalized
intermediate Jacobian
\begin{equation} J^n(X) = \rm
H^{2n-1}(X_{\mathrm{an}},{\mathbb C}){\Large /}\rm
H^{2n-1}(X_{\mathrm{an}},{\mathbb Z}(n)) + \mathrm{F}^n\rm
H^{2n-1}(X_{\mathrm{an}},{\mathbb C}),\end{equation}  but by the
Jacobians of the curves that are the building blocks of the
middle-dimensional cohomology $\rm H^{\mathrm{dim}_{{\mathbb C}}
X}(X)$. These Jacobians themselves do not admit complex
multiplication, unlike the situation in the elliptic case, but
instead split into different factors which admit different types
of complex multiplication, in general. Furthermore the ring class
field can be generalized to be the field of moduli, and we can
consider also points on the abelian variety that are of finite
order, i.e. torsion points, and the field extensions they
generate.

This allows us to answer a question posed in \cite{m98} which
asked whether the absolute Galois group $\mathrm{Gal}(\bar{K}/K)$
could play a role in the context of $N=2$ compactifications of
type IIB strings. This is indeed the case. Suppose we have given
an abelian variety $A$ defined over a field $K$ with complex
multiplication by a field $E$. Then there is an action of the
absolute Galois group $\mathrm{Gal}(\bar{K}/K)$ of the closure
$\bar{K}$ of $K$ on the torsion points of $A$. This action is
described by a Hecke character which is associated to the fields
$(K,E)$ \cite{s71}.

We have mentioned already that in general the (Griffiths)
intermediate Jacobian is only a torus, not an abelian variety.
Even in those cases it is however possible to envision the
existence of motives via abelian varieties associated to a variety
$X$. Consider the Chow groups $\mathrm{CH}^p(X)$ of codimension
$p$ cycles modulo rational equivalence and denote by
$\mathrm{CH}^p(X)_{\mathrm{Hom}}$ the subgroup of cycles
homologically equivalent to zero. Then there is a homomorphism,
called the Abel-Jacobi homomorphism, which embeds
$\mathrm{CH}^p(X)_{\mathrm{hom}}$ into the intermediate Jacobian
\begin{equation} \phi: \mathrm{CH}^p(X)_{\mathrm{hom}} \longrightarrow
J^p(X).\end{equation} The image of $\phi$ on the subgroup
${\mathcal A}^p(X)$ defined by cycles algebraically equivalent to
zero does in fact define an abelian variety, even if $J^p(X)$ is
not an abelian variety but only a torus \cite{s79}. Hence we can
ask whether attractor varieties are distinguished by Abel-Jacobi
images which admit complex multiplication.

Even more general, we can formulate this question in the framework
of motives because of Delign's conjecture. Thinking of motives as
universal cohomology theories, it is conceivable that attractor
varieties lead to motives in the abelian category with (potential)
complex multiplication. The standard cycle class map construction
of $\mathrm{CH}^p(X)_{\mathrm{hom}}$ is replaced by the first term
of a (conjectured) filtration in the resulting K-theory.

Putting everything together we see that the two separate
discussions in \cite{m98} characterizing toroidal attractor
varieties on the one hand, and Calabi-Yau hypersurfaces on the
other, are just two aspects of our way of looking at this problem.
This is the case precisely because of Deligne's period conjecture
which relates the field of the periods to the field of complex
multiplication via the L-function of the variety (or motive). Thus
a very pretty unified picture emerges.

\vskip .2truein

{\bf Acknowledgements}. M.L. and R.S. were supported in part by
the NSF under Grant No. PHY99-07949. They are grateful to the
Kavli Institute for Theoretical Physics, Santa Barbara, for
hospitality and support through KITP Scholarships during the
course of part of this work.
 Their work has also been supported in part by NATO under Grant
No. CRG 9710045. The work of V.P. was supported in part by the NSF
under Grant No. PHY98-02484.

\end{document}